\documentclass[9pt,amstex]{article}
\usepackage{amssymb,amsmath} 
\usepackage{latexsym}
\usepackage[usenames,dvipsnames]{color}

\newtheorem{dfn}{Definition}[section] 
\newtheorem{rmk}{Remark}[section]

\newtheorem{prop}{Proposition}[section]

\def\cyclic{\mathop{\kern0.9ex{{+}
\kern-2.2ex\raise-.28ex\hbox{\Large\hbox
{$\circlearrowright$}}}}\limits}

\def\buildrel#1_#2^#3{\mathrel{\mathop{\kern 0pt#1}\limits_{#2}^{#3}}}

\newcommand{\Pf}{{\em Proof}. }
\newcommand{\EPf}
{%
\mbox{}%
\nolinebreak%
\hfill%
\rule{2mm}{2mm}%
\medbreak%
\par%
}

\newcommand{\End}{\mbox{$\mathtt{End}$}}

\newcommand{\id}{\mbox{$\mathtt{id}$}}

\newcommand{\Ad}{\mbox{$\mathtt{Ad}$}}

\newcommand{\C}{\mathbb C} 
\newcommand{\E}{\mathbb E}

\newcommand{\bC}{{\bf C}} 
\newcommand{\bB}{{\bf B}}

\newcommand{\bH}{{\bf H}}

\newcommand{\T}{\mathbb T}

\newcommand{\F}{\mathbb F}

\newcommand{\R}{\mathbb R}

\newcommand{\g}{{\mathfrak{g}}{}}

\renewcommand{\k}{{\mathfrak{k}}{}}

\newcommand{\fL}{{\mathfrak{L}}{}}

\newcommand{\q}{{\mathfrak{q}}{}}

\newcommand{\ddto}{{\frac{{\rm d}}{{\rm d}t}|_{0}}{}}

\newcommand{\fu}{{\mathfrak{u}}{}} 
\newcommand{\h}{{\mathfrak{h}}{}} 
 
\renewcommand{\b}{{\mathfrak{b}}{}} 

\newcommand{\CO}{{\cal O}{}}

\newcommand{\CE}{\mathcal E}

\newcommand{\CL}{{\cal L}{}}

\newcommand{\CH}{\mathcal H}
\newcommand{\CB}{\mathcal B}

\newcommand{\CF}{\mathcal F}

\newcommand{\ubf}{U_{\mbox{\rm \tiny{BF}}}}
\newcommand{\ukw}{U_{\mbox{\rm \tiny{KW}}}}

\title{Bargmann-Fock realization of the noncommutative torus}
\author{Nestor Anzola Kibamba and Pierre Bieliavsky\\ IRMP, UCLouvain, Belgium}
\date{}

\begin{document}
\maketitle
\begin{abstract}
 We give an interpretation of the Bargman transform as a correspondence
between state spaces that is analogous to  commonly considered intertwiners in representation theory of finite groups.
We observe that the non-commutative torus is nothing else that the range of the star-exponential for the Heisenberg group
within the Kirillov's orbit method context. We deduce from this a realization of the non-commutative torus as acting 
on a Fock space of entire functions.
\end{abstract}

\section{Introduction}
For systems of finite degrees of freedom, there are two main quantization procedures, the other ones being variants of them.
The first one consisting in pseudo-differential calculus (see e.g. \cite{Ho}). The second one relies on geometric quantization (see e.g. \cite{Wo}). The main difference
between the two lies in the types of polarizations on which they are based. Pseudo-differential calculus
is based on the existence of a ``real polarization", while geometric quantization often uses ``complex polarizations".
There are no systematic ways to compare the two. Although in some specific situations, this comparison is possible.
This is what is investigated in the present work.

The aim of the small note is threefold. First we give an interpretation of the Bargman transform as a correspondence
between state spaces that is analogous to  commonly considered intertwiners in representation theory of finite groups.
Second, we observe that the non-commutative torus is nothing else that the range of the star-exponential for the Heisenberg group
within the Kirillov's orbit method context. Third, we deduce from this a realization of the non-commutative torus as acting 
on a Fock space of entire functions. The latter relates the classical approach to the non-commutative torus in the context 
of Weyl quantization to its realization, frequent in the Physics literature, in terms of the canonical quantization.

\vspace{2mm}

\noindent{\bf Acknowledgment} The authors thank the Belgian Scientific Policy (BELSPO) for support under IAP grant
P7/18 DYGEST.

\section{Remarks on the geometric quantization of co-adjoint orbits} We let $\CO$ be a co-adjoint orbit of a connected Lie group $G$ in the dual $\g^\star$ of its Lie algebra $\g$.
We fix a base point $\xi_o$ in $\CO$ and denote by $K:=:G_{\xi_o}$ its stabilizer in $G$ (w.r.t. the co-adjoint action). 
We assume $K$ to be connected.
Denote by $\k$ the Lie sub-algebra of $K$ in $\g$. Consider the $\R$-linear map
$\xi_o:\k\to\fu(1)=\R:Z\mapsto<\xi_o,Z>\;.$
Since the two-form $\delta\xi_o:\g\times\g\to\R:(X,Y)\mapsto<\xi_o,[X,Y]>$ identically vanishes on $\k\times\k$,
the above mapping is a character of $\k$.
Assume the above character exponentiates to $K$ as a unitary character (Kostant's condition):
$\chi:K\to U(1)\quad(\chi_{\star e}=i\xi_o)\;.$
One then has an action of $K$ on $U(1)$ by group automorphisms:
$K\times U(1)\to U(1):(k,z)\mapsto\chi(k)z\;.$
The associated circle bundle $Y\;:=\;G\times_KU(1)$ is then naturally a $U(1)$-principal bundle over the orbit $\CO$:
$\pi:Y\to\CO:[g,z]\mapsto\Ad^\flat_g(\xi_o)$
(indeed, one has the well-defined $U(1)$-right action $[g,z].z_0:=[g,zz_0]$).

\noindent The data of the character yields a connection one-form $\varpi$ in $Y$. Indeed,
the following formula defines a left-action of $G$ on $Y$:
$g_0.[g,z]\;:=\;[g_0g,z]\;.$

\noindent When $\xi_o|_\k$ is non-trivial, the latter is transitive: $\bC_g(k).[g,z]=[gk,z]=[g,\chi(k)z]$ and $\pi(g_0.[g,z])=g_0.\pi(g)$. We then set
\begin{equation}\label{1FORM}
\varpi_{[g,z]}(X^\star_{[g,z]})\;:=\;-\,<\Ad^\flat_g\xi_o,X>
\end{equation}
with, for every $X\in\g$: $X^\star_{[g,z]}\;:=\;\ddto\exp(-tX).[g,z]$. The above formula (\ref{1FORM})  defines a 1-form. Indeed, an element
$X\in\g$ is such that $X^\star_{[g,z]}=0$ if and only if $\Ad_{g^{-1}}X\in\ker(\xi_0)\cap\k$. It is a connection form because for every $z_0\in U(1)$, one has
$(z_0^\star\varpi)_{[g,z]}(X^\star)=\varpi_{[g,zz_0]}(z_{0\star_{[g,z]}}X^\star)=\varpi_{[g,zz_0]}(X^\star)=-<\Ad^\flat_g\xi_o,X>=\varpi_{[g,z]}(X^\star)$.
At last, denoting by $\iota_y\;:=\;\ddto y.e^{it}$ the infinitesimal generator of the circle action on $Y$, one has
$\iota_{[g,z]}\;=\;-(\Ad_gE)^\star_{[g,z]}$
where $E\in\k$ is such that $<\xi_0,E>=1$. Indeed, 
$\ddto [g,z].e^{it}=\ddto [g,e^{it}z]=\ddto [g\exp(tE),z]=\ddto [\exp(t\Ad_gE)g,z]\;.$
Therefore: $\varpi_{[g,z]}(\iota)=<\Ad_g^\flat\xi_0,\Ad_gE>=1$.

\noindent The curvature $\Omega^{\varpi}:={\rm d}\varpi+[\varpi,\varpi]={\rm d}\varpi$ of that connection equals the lift $\pi^\star(\omega^\CO)$ to $Y$
of the KKS-symplectic structure $\omega^\CO$ on $\CO$ because $X^\star_{[g,z]}.\varpi(Y^\star)=<\Ad_g^\flat\xi_o,[X,Y]>$.

\subsection{Real polarizations} We now relate Kirillov's polarizations to Souriau's Planck condition. We consider a {\bf partial polarization} affiliated to $\xi_o$ i.e. a sub-algebra $\b$ of $\g$ that is normalized by $\k$ and  maximal for the property of being isotropic w.r.t. the two-form $\delta\xi_o$ on $\g$ defined as
$\delta\xi_o(X,Y)\;:=\;<\xi_o,[X,Y]>$. We assume that the analytic (i.e. connected) Lie sub-group $\bB$ of $G$ whose Lie algebra is $\b$ is closed. We denote by $Q\;:=\;G/\bB$ the corresponding quotient manifold. Note that one necessarily has $K\subset\bB$, hence the fibration
$p:\CO\to Q:\Ad^\flat_g(\xi_o)\mapsto g\bB$. The distribution $\fL$ in $T\CO$ tangent to the fibers is isotropic w.r.t. the KKS form. Its $\varpi$-horizontal lift $\overline{\fL}$ in $T(Y)$, being integrable, induces a {\bf Planck foliation} of $Y$ (cf. p. 337 in Souriau's book \cite{S}).

\noindent Usually, Kirllov's representation space consists in a space of sections of the associated complex line bundle $\E_\chi:=G\times_\chi\C\to Q$ where $\chi$ is viewed as a unitary character of $\bB$. While Kostant-Souriau representation space rather 
consists in a space of sections of the line bundle $Y\times_{U(1)}\C\to\CO$. One therefore looks for a morphism between these spaces.
To this end we first observe that the circle bundle $G\times_\chi U(1)\to Q$ is a principal $U(1)$-bundle similarly as  $Y$ is over $\CO$. Second, we note the complex line bundle isomorphism over $Q$:
\begin{equation}\label{U(1)BUNDLEISO}
(G\times_\chi U(1))\times_{U(1)}\C\to G\times_\chi\C:[[g_0,z_0],z]\mapsto[g,z_0z]\;.
\end{equation}
Third, we have the morphism of $U(1)$-bundle:
$$
\begin{array}{ccc}
Y&\stackrel{\tilde{p}}{\longrightarrow}&G\times_\chi U(1)\\
\downarrow&&\downarrow\\
\CO&\stackrel{{p}}{\longrightarrow}&Q
\end{array}
$$
where $\tilde{p}([gk,\chi(k^{-1})z_0])\;:=\;[gb,\chi(b^{-1})z_0]$.
This induces a linear map between equivariant functions:
$\tilde{p}^\star: C^\infty(G\times_\chi U(1),\C)^{U(1)}\to C^\infty(Y,\C)^{U(1)}$
which, red through the isomorphism (\ref{U(1)BUNDLEISO}), yields a natural $G$-equivariant linear embedding
$\Gamma^\infty(Q,\E_\chi)\to\Gamma^\infty(\CO,Y\times_{U(1)}\C)$
whose image coincides with the Planck space.

\subsection{Complex polarizations} We first note that for every element $X\in\g$, the $\varpi$-horizontal lift of $X^\star_\xi$ at $y=[g,z]\in Y$
is given by
$h_y(X^\star_\xi)\;:=\;X^\star_y+<\xi,X>\iota_y\;.$
Indeed, $\pi_{\star y}(X^\star_y)=X^\star_\xi$ and $\varpi_y(X^\star_y+<\xi,X>\iota_y)=-<\xi,X>+<\xi,X>=0$.

\noindent Therefore, for every smooth section $\varphi$ of the associated complex line bundle $\F\;:=\;Y\times_{U(1)}\C\to\CO$, 
denoting by $\nabla$ the covariant derivative in $\F$ associated to $\varpi$ and by $\hat{\varphi}$ the $U(1)$-equivariant function on $Y$ representing $\varphi$, one has
$\widehat{\nabla_{X^\star}\varphi}(y)=X^\star_y\hat{\varphi}\,-\,i\,<\xi,X>\hat{\varphi}(y)\;.$
Now let us assume that our orbit $\CO$ is \emph{pseudo-Kahler} in the sense that it is equipped with a $G$-invariant
$\omega^\CO$-compatible (i.e. $J\in\mbox{\rm Sp}(\omega^\CO)$) almost complex structure $J$. Let us denote by 
$T_\xi(\CO)^\C\;=\;T_\xi^{0}(\CO)\,\oplus\,T^1_\xi(\CO)$
the $(-1)^{0,1} i$-eigenspace decomposition of the complexified tangent space $T_\xi(\CO)$ w.r.t. $J_\xi$.
One observes the following descriptions
$T_\xi^{0,1}(\CO)\;=\;>X+(-1)^{0,1}iJX<_{X\in T_\xi(\CO)}$
where the linear span is taken over the complex numbers. We note also that $\dim_\C T_\xi^{0,1}(\CO)=\frac{1}{2}\dim_\R\CO$.

\noindent In that context, a smooth section $\varphi$ of $\F$ is called {\bf polarized} when $\nabla_Z\varphi=0$
for every $Z\in T^0(\CO)$. The set $\mbox{\rm hol}(\F)$ of polarized sections is a complex sub-space of 
$\Gamma^\infty(\F)$. Moreover, it carries a natural linear action of $G$. Indeed, the group $G$ acts on 
$\Gamma^\infty(\F)$ via $\widehat{U(g)\varphi}:=(g^{-1})^\star\hat{\varphi}$. The fact that both $\varphi$ and $J$
are $G$-invariant then  implies that $\mbox{\rm hol}(\F)$ is a $U$-invariant sub-space of $\Gamma^\infty(\F)$. The linear representation
$U:G\to\End(\mbox{\rm hol}(\F))$
is called the {\bf Bargman-Fock representation}.

\section{The Heisenberg group} We consider a symplectic vector space $(V,\Omega)$ of dimension $2n$ and its associated 
Heisenberg Lie algebra: $\h_n:=V\oplus\R E$ whose Lie bracket is given by $[v,v']:=\Omega(v,v')E$ for all $v,v'\in V$, the element $E$
being central. The corresponding connected simply connected Lie group $\bH_n$ is modeled on $\h_n$ with group
law given by $g.g':=g+g'+\frac{1}{2}[g.g']$. 

\noindent Within this setting, one observes that the exponential mapping is the identity map on $\h_n$. The symplectic structure defines an isomorphism ${}^\flat:V\to V^\star$ by
${}^\flat v(v'):=\Omega(v,v')$. The latter extends to a linear isomorphism ${}^\flat:\h_n\to\h_n^\star$ where we set
${}^\flat E(v+zE):=z$. 

\noindent Now one has \begin{equation}\label{ADFLAT}\Ad^\flat_{v+zE}({}^\flat v_0+\mu{}^\flat E)={}^\flat(v_0-\mu v)+\mu{}^\flat E\;.\end{equation}
Indeed, in the case of the Heisenberg group, the exponential mapping coincides with the identity mapping. Hence, for every
$v_1+z_1E\in\h_n$:
\begin{eqnarray*}
&& <\Ad^\flat_{v+zE}({}^\flat v_0+\mu{}^\flat E),v_1+z_1E>=<{}^\flat v_0+\mu{}^\flat E,\Ad_{-v-zE}(v_1+z_1E)>\\
&=&\ddto<{}^\flat v_0+\mu{}^\flat E,(-v-zE).(tv_1+tz_1E).(v+zE)>\;,
\end{eqnarray*}
which in view of the above expression for the group law immediately yields (\ref{ADFLAT}).

\noindent Generic orbits $\CO$ are therefore affine hyperplanes parametrized by $\mu\in\R_0$. Setting $\xi_0:=\mu{}^\flat E$, every real
polarization $\b$ corresponds to a choice of a Lagrangian sub-space $\fL$ in $V$: $\b=\fL\oplus\R E$. Note in particular, that 
the polarization is an ideal in $\h_n$. Choosing a Lagrangian sub-space
$\q$ in duality with $\fL$ in $V$ determines an Abelian sub-group $Q=\exp(\q)$ in $\bH_n$ which splits the exact sequence
$\bB\to\bH_n\to\bH_n/\bB=:Q$ i.e. $\bH_n=Q.\bB$. The stabilizers all coincide (in the generic case) with the center $K:=\R E$
of $\bH_n$ and one has the global trivialization $\CO\to\bH_n:{}^\flat v+\xi_0\mapsto v$. 

\subsection{Representations from real and complex polarizations}
This yields the linear isomorphism
$\Gamma^\infty(\E_\chi)\to C^\infty(Q):u\mapsto\hat{u}|_Q=:\tilde{u}$ under which the $\bH_n$-action reads
$U_{\mbox{\rm \tiny{KW}}}(qb)\tilde{u}(q_0)=e^{i\mu(z+\Omega(q-q_0,p))}\tilde{u}(q_0-q)$ with $b=p+zE\,,\,p\in\exp(\fL)$. The latter induces a unitary 
representation on $L^2(Q)$.

\noindent The isomorphism $\C^n=\q\oplus i\fL\to V:Z=q+ip\mapsto q+p$ yields an $\bH_n$-equivariant global complex
coordinate system on the orbit $\CO$ through $Z\mapsto\Ad^\flat_{q+p}\xi_0$. The map $V\times U(1)\to Y:(v_0,z_0)\mapsto[v_0,z_0]$
consists in a global trivialization of the bundle $Y\to\CO$ under the isomorphism $V\simeq\CO$ described above. Hence the linear isomorphism $\Gamma^\infty(\F)\to C^\infty(V):\varphi\mapsto\tilde{\varphi}:=\hat{\varphi}(\,.\,,\,1)$. At the level of the trivialization the left-action
of $\bH_n$ on $Y$ reads: $g.(v_0,z_0)=(v_0+q+p,e^{i\mu(z+\frac{1}{2}\Omega(q,p)+\frac{1}{2}\Omega(q+p,v_0))}z_0)$. Also, 
the representation is given by $g.\tilde{\varphi}(v_0)=e^{i\mu(z+\frac{1}{2}\Omega(q,p)+\frac{1}{2}\Omega(q+p,v_0))}\tilde{\varphi}(v_0-p-q)$.
Choosing a basis $\{f_j\}$ of $\q$ and setting $\{e_j\}$ for the corresponding dual basis of $\fL$, one has
$\partial_{\overline{Z}^j}=-\frac{1}{2}(f^\star_j+ie^\star_j)$. Within the trivialization, the connection form corresponds to $\varpi=\frac{\mu}{2}(p^j{\rm d}q_j-q^j{\rm d}p_j)+\iota_\star$. A simple computation then yields:
$\mbox{\rm hol}(\F)\;=\;\{\varphi_f\;:\;\C^n\to\C:z\mapsto e^{-\frac{\mu}{4}|z|^2}\,f(z)\;\quad (\,f\;\mbox{\rm entire}\,)\;\}\;.$
Note that provided $\mu>0$, the space $\mbox{\rm hol}(\F)$ naturally contains the pre-Hilbert space:
$\CL^2_{\mbox{\rm hol}}(\F)\;:=\;\{\varphi_f\;:\;<\varphi_f,\varphi_f>\;:=\;\int_{\C^n}e^{-\frac{\mu}{2}|z|^2}\,\left| f(z)\right|^2\,{\rm d}q\,{\rm d}p\,<\,\infty\;\}\;.$
The above sub-space turns out to be invariant under the representation $U$ of $H_n$. The latter is seen to be
unitary and irreducible. For negative $\mu$, one gets a unitary representation by considering anti-polarized sections corresponding to 
anti-holomorphic functions. Note that within complex coordinates, the action reads $U_{\mbox{\rm \tiny{BF}}}(g)\tilde{\varphi}(Z_0):=g.\tilde{\varphi}(Z_0)=e^{i\mu(z+\frac{1}{2}\mbox{\rm Im}(\frac{1}{2}{Z}^2+\overline{Z}Z_0))}\tilde{\varphi}(Z_0-Z)$.

\subsection{Intertwiners and the Bargmann transform}
\noindent We know (cf. Stone-von Neumann's theorem) that in the case of the Heisenberg group, representations constructed either
via complex or via real polarizations are equivalent. In order to exhibit intertwiners, we make the following observation.
\begin{prop}
Let $G$ be a Lie group with left-invariant Haar measure ${\rm d}g$ and $(\CH,\rho)$ and $(\CH',\rho')$ be square-integrable unitary representations. Assume furthermore the continuity of the associated bilinear forms $\CH\times\CH\to L^2(G):(\varphi_1,\varphi_2)\to[g\mapsto<\varphi_1|\rho(g)\varphi_2>]$ and $\CH'\times\CH'\to L^2(G):(\varphi'_1,\varphi'_2)\to[g\mapsto<\varphi_1'|\rho'(g)\varphi'_2>]$. Fix ``mother states" $|\eta>\in\CH$ and $|\eta'>\in\CH'$.
For every element $g\in G$, set $|\eta_g>:=\rho(g)|\eta>$ and $|\eta'_g>:=\rho'(g)|\eta'>$.
Then the following formula 
\begin{equation}\label{INTERTWINER}
T:=\int_G\,|\eta'_g><\eta_g|\,{\rm d}g
\end{equation}
formally defines an intertwiner from $(\CH,\rho)$ to $(\CH',\rho')$. 
\end{prop}
\Pf
First let us observe that square-integrability and Cauchy-Schwartz inequality on $L^2(G)$  imply that for all $\varphi\in\CH$ and $\varphi'\in\CH'$ the element $[g\mapsto<\varphi'|\eta'_g><\eta_g|\varphi>]$ is well defined as an element of $L^1(G)$. Moreover the continuity of the above mentioned bilinear forms insures the continuity  of the bilinear map $\CH\times\CH'\to\C$ defined by integrating the later. This is in this sense that we understand formula (\ref{INTERTWINER}). Now for all $|\varphi>\in\CH$ and $g_0\in G$, one has
$T|\varphi_{g_0}>=\int|\eta'_g><\eta_g|\varphi_{g_0}>\,{\rm d}g=\int|\eta'_g><\eta_{g_0^{-1}g}|\varphi>\,{\rm d}g=
\int|\eta'_{g_0g}><\eta_{g}|\varphi>\,{\rm d}g=\int\rho'(g_0)|\eta'_{g}><\eta_{g}|\varphi>\,{\rm d}g=\rho'(g_0)T|\varphi>$.
\EPf

\noindent In our present context of the Heisenberg group, the integration over $G$ should rather be replaced by an integration
over the orbit (which does not correspond to a sub-group of $\bH_n$). But the above argument essentially holds the same way up to a slight modification by a pure phase. Namely, 
\begin{prop}
Fix $\tilde{\varphi}^0\in\CL^2_{\mbox{\rm hol}}(\F)$ and $\tilde{u}^0\in L^2(Q)$. 
For every $v\in V$, set $\tilde{\varphi}^0_v:=\ubf(v)\tilde{\varphi}^0$ and $\tilde{u}^0_v:=\ukw(v)\tilde{u}^0$. Then, setting
$$
T_{\varphi^0,u^0}\;:=\;\int_V\,|\tilde{\varphi}^0_v><\tilde{u}^0_v|\,{\rm d}v
$$
formally defines a $V$-intertwiner between $\ukw$ and $\ubf$.
\end{prop}
\Pf Observe firstt that for every 
$w\in V$, one has $\tilde{\varphi}^0_{wv}=e^{\frac{i\mu}{2}\Omega(w,v)}\tilde{\varphi}^0_{w+v}$ and similarly for $\tilde{u}^0$. Also for every $\tilde{u}\in L^2(Q)$, one has
\begin{eqnarray*}
T_{\varphi^0,u^0}\ukw(w)\tilde{u}&=&\int_V\,|\tilde{\varphi}^0_v><\tilde{u}^0_v|\tilde{u}_w>\,{\rm d}v=
\int_V\,|\tilde{\varphi}^0_v><\tilde{u}^0_{w^{-1}v}|\tilde{u}>\,{\rm d}v\\
&=&
\int_V\,|\tilde{\varphi}^0_v><e^{-\frac{i\mu}{2}\Omega(w,v)}\tilde{u}^0_{v-w}|\tilde{u}>\,{\rm d}v=
\int_V\,e^{\frac{i\mu}{2}\Omega(w,v)}|\tilde{\varphi}^0_v><\tilde{u}^0_{v-w}|\tilde{u}>\,{\rm d}v\\
&=&
\int_V\,e^{\frac{i\mu}{2}\Omega(w,v+w)}|\tilde{\varphi}^0_{v+w}><\tilde{u}^0_{v}|\tilde{u}>\,{\rm d}v=
\int_V\,e^{\frac{i\mu}{2}\Omega(w,v)}|\tilde{\varphi}^0_{v+w}><\tilde{u}^0_{v}|\tilde{u}>\,{\rm d}v\\
&=&
\int_V\,|\tilde{\varphi}^0_{wv}><\tilde{u}^0_{v}|\tilde{u}>\,{\rm d}v=\ubf(w)T_{\varphi^0,u^0}\tilde{u}\;.
\end{eqnarray*}
Now, one needs to check whether the above definition makes actual sense. The special choices $\tilde{\varphi}^0(Z):=\varphi_1(Z)=e^{-\frac{\mu}{4}|Z|^2}$ and
$\tilde{u}^0(q):=e^{-\alpha q^2}$ reproduce the usual Bargman transform. Indeed, first observe that
$<\tilde{u}^0_v,\tilde{u}>=e^{\frac{-i\mu}{2}qp}\int_Qe^{i\mu q_0p-\alpha(q_0-q)^2}\tilde{u}(q_0){\rm d}q_0$ and 
$\tilde{\varphi}^0_v(v_1)=e^{i\frac{\mu}{2}(qp_1-q_1p)-\frac{\mu}{4}((q_1-q)^2+(p_1-p)^2)}$. This leads to
$T\tilde{u}(v_1)=\int{\rm d}q_0{\rm d}q{\rm d}p \;e^{\frac{i\mu}{2}((p-p_1)(2q_0-q_1-q)+(2q_0-q_1)p_1)}e^{-\alpha(q_0-q)^2-\frac{\mu}{4}((q_1-q)^2+(p_1-p)^2)}\tilde{u}(q_0)$. From the fact that $\int e^{-ixy}e^{-\frac{x^2}{2\sigma^2}}{\rm d}x=e^{-\frac{\sigma^2}{2}y^2}$, we get: $$T\tilde{u}(v_1)=\left(2\sqrt{\frac{\pi}{\mu}}\right)^n\int{\rm d}q_0{\rm d}q \;e^{\frac{i\mu}{2}(2q_0-q_1)p_1}e^{-\alpha(q_0-q)^2-\frac{\mu}{4}(q_1-q)^2}e^{-\frac{\mu}{4}(2q_0-q_1-q)^2}\tilde{u}(q_0)\;.$$
The formula $\int e^{-\frac{2}{\sigma^2}(q-q_0)^2}{\rm d}q=(\frac{\sqrt{2\pi}}{\sigma})^n$ yields
$$T\tilde{u}(v_1)=\left(2\pi\sqrt{\frac{\alpha+\frac{\mu}{4}}{\mu}}\right)^n\int{\rm d}q_0 \;e^{\frac{i\mu}{2}(2q_0-q_1)p_1}e^{-\frac{\mu}{2}((q_1-q_0)^2+\frac{1}{2}q_0^2)}\tilde{u}(q_0)\;.$$
Setting $Z_1:=q_1+ip_1$ leads to
$$T\tilde{u}(v_1)=e^{-\frac{\mu}{4}|Z_1|^2}\left(2\pi\sqrt{\frac{\alpha+\frac{\mu}{4}}{\mu}}\right)^n\int{\rm d}q_0 \;e^{-\frac{\mu}{4}(Z_1-q_0)(Z_1-3q_0)}\tilde{u}(q_0)\;.$$
\EPf
\begin{rmk}The usual Bargman transform (see Folland \cite{F} page 40) equals the latter when $\frac{\alpha}{3}=\frac{\pi}{2}=\frac{\mu}{4}$.\end{rmk}

\section{Star-exponentials and noncommutative tori}
\subsection{Recalls on Weyl calculus}
It is well known that the Weyl-Moyal algebra can be seen as a by-product of the Kirillov-Weyl representation. In \cite{BGT},
this fact is realized in terms of the natural symmetric space structure on the coadjoint orbits of the Heisenberg group.
This construction is based on the fact that the Euclidean centered symmetries on $V=\R^{2n}\simeq\CO$ naturally ``quantize as phase
symmetries". More precisely, at the level of the Heisenberg group the flat symmetric space structure on $V$ is encoded 
by the involutive automorphism $$\sigma:\bH_n\to\bH_n:v+zE\mapsto-v+zE\;.$$ The latter yields an involution
of the equivariant function space $$\sigma^\star:C^\infty(\bH_n,\C)^B\to C^\infty(\bH_n,\C)^B$$ which induces by restriction to
$Q$ the unitary phase involution: $$\Sigma:L^2(Q)\to L^2(Q):\tilde{u}\mapsto[q\mapsto\tilde{u}(-q)]\;.$$ Observing that the associated map
$$\bH_n\to U(L^2(Q)):g\mapsto U_{\mbox{\rm \tiny{KW}}}(g)\,\Sigma\,U_{\mbox{\rm \tiny{KW}}}(g^{-1})$$ is constant on the 
lateral classes of the stabilizer group $\R E$, one gets an $\bH_n$-equivariant  mapping: $$\Xi_\mu:V\simeq\bH_n/\R E\to U(L^2(Q)): 
v\mapsto U_{\mbox{\rm \tiny{KW}}}(v)\,\Sigma\,U_{\mbox{\rm \tiny{KW}}}(v^{-1})$$ which at the level of functions yields the
so-called ``Weyl correspondance" valued in the $C^\star$-algebra of bounded operators on $L^2(Q)$:
$$
\Xi_\mu:L^1(V)\to\CB(L^2(Q)):F\mapsto\int_VF(v)\,\Xi_\mu(v)\,{\rm d}v\;.
$$
The above mapping uniquely extends from the space of compactly supported functions $C^\infty_c(V)$ to an injective  continuous
map defined on the Laurent Schwartz B-space $\CB(V)$ of complex valued smooth functions on $V$ whose partial derivatives at every order are bounded:
$$
\Xi_\mu:\CB(V)\to\CB(L^2(Q))\;,
$$
expressing in particular that the quantum operators associated to classical observables in the B-space are $L^2(Q)$-bounded
in accordance with the classical Calder\`on-Vaillancourt theorem. 

\noindent It turns out that the range of the above map is a sub-algebra of $\CB(L^2(Q))$. The Weyl-product $\star_\mu$ on $\CB(V)$ is then defined by structure transportation:
$$
F_1\star_\mu F_2\;:=\;\Xi_\mu^{-1}\left(\Xi_\mu(F_1)\circ\Xi_\mu(F_2)\right)
$$
whose asymptotics in terms of powers of $\theta\;:=\;\frac{1}{\mu}$ consists in the usual formal Moyal-star-product:
$$
F_1\star_\mu F_2\;\sim\;F_1\star^0_\theta F_2\;:=\;\sum_{k=0}^\infty\frac{1}{k!}\left(\frac{\theta}{2i}\right)^k\Omega^{i_1j_1}...
\Omega^{i_kj_k}\partial^k_{i_1...i_k}F_1\,\partial^k{j_1...j_k}F_2\;.
$$
The resulting associative 
algebra $(\CB(V),\star_\mu)$ is then both Fr\'echet (w.r.t. the natural Fr\'echet topology on $\CB(V)$) and pre-$C^\star$ 
(by transporting the operator norm from $\CB(L^2(Q))$).

\subsection{Star-exponentials}
Combining the above mentioned results of \cite{BGT} and well-known results on star-exponentials (see e.g. \cite{A}), 
we observe that the heuristic consideration of the series:
$$
\mbox{Exp}_{\theta}(F)\;:=\;\sum_{k=0}^\infty\frac{1}{k!}\left(\frac{i}{\theta}F\right)^{\star_\mu k}
$$
yields a well-defined group homomorphism:
$$
\CE_\theta:\bH_n\to(\CB(V),\star_\mu):g\mapsto\mbox{Exp}_{\theta}(\lambda_{\log g})
$$
where $\lambda$ denotes the classical linear moment:
$$
\h_{n}\to C^\infty(V):X\mapsto[v\mapsto<\Ad_v^\flat\xi_0,X>]\;.
$$
Indeed, in this case, if $F$ depends only either on the $q$-variable or on the $p$-variable then 
the above star-exponential just coincides with the usual exponential:
$\mbox{Exp}_{\theta}(F)=\exp\left(\frac{i}{\theta}F\right)$. In particular, for $x$ either in $Q$ or $\exp\fL$ we find:
$$
\left(\CE_\theta(x)\right)(v)\;=\;e^{\frac{i}{\theta^2}\Omega(x,v)}\;;
$$
while for $x=zE$, we find the constant function:
$$
\left(\CE_\theta(zE)\right)(v)\;=\;e^{\frac{zi}{\theta^2}}\;.
$$
From which we conclude that $\CE_\theta$ is indeed valued in $\CB(V)$.

\subsection{An approach to the non-commutative torus}
For simplicity, restrict to the case $n=1$ and consider $\Omega$-dual basis elements $e_q$ of $\q$ and $e_p$  of $\fL$.
From what precedes we observe that the group elements
$$
U\;:=\;\CE_\theta(\exp(\theta^2e_q))\;=\;e^{ip}\quad\mbox{ and }\quad V\;:=\;\CE_\theta(\exp(\theta^2e_p))\;=\;e^{-iq}
$$
behave inside the image group $\CE_\theta(\bH_1)\;\subset\;\CB(\R^2)$ as
$$
U\,V\;=\;e^{i\theta}V\,U
$$
(where we omitted to write $\star_\mu$).

\noindent In particular, we may make the following 
\begin{dfn} Endowing $(\CB(\R^2),\star_\mu)$ with its pre-$C^\star$-algebra structure (coming from $\Xi_\mu$),
the complex linear span  inside $\CB(\R^2)$ of the sub-group of $\CE_\theta(\bH_1)$
generated by elements $U$ and $V$ underlies a pre-$C^\star$-algebra $\T^\circ_\theta$ that completes as the non-commutative 2-torus
$\T_\theta$.
\end{dfn}

\subsection{Bargman-Fock realization of the non-commutative torus}
Identifying elements of $\R^2=V\subset\bH_1$ with complex numbers as before, we compute that
$$
T\circ\Sigma\;=\;-\id^\star\circ T\quad\mbox{and }\quad \mbox{\rm BF}_\mu(Z)\tilde{\varphi}(Z_0)\;:=\;T\Xi_\mu(Z)T^{-1}\tilde{\varphi}(Z_0)\;=\;e^{i\mu\,\mbox{\rm Im}(\overline{Z}Z_0)}
\tilde{\varphi}(2Z-Z_0)\;.
$$
By structure transportation, we define the following correspondence:
$$
 \mbox{\rm BF}_\mu:\CB(\C)\to\CB(L^2_{\mbox{\rm hol}}(\F))
$$
as the unique continuous linear extension of 
$$
C^\infty_c(\C)\to\CB(L^2_{\mbox{\rm hol}}(\F)):F\mapsto\int_\C F(Z)\,\mbox{\rm BF}_\mu(Z)\,{\rm d}Z\;.
$$
\begin{prop}
Applied to an element  in $\CE_\theta(\bH_1)$ of the form $F_X(Z)\;=\;e^{i\alpha \, \mbox{\rm Im}(\overline{Z}X)}$ with $X\in\C$ and $\alpha\in\R$, one has:
$$
\mbox{\rm BF}_\mu(F_X)
\tilde{\varphi}(Z_0)\;=\;e^{\frac{i\alpha}{2}\mbox{\rm Im}(\overline{Z_0}X)}\tilde{\varphi}(\mu Z_0+\alpha X)\;.
$$
\end{prop}
\Pf
A small computation leads to the formula:
$$
 \mbox{\rm BF}_\mu(F)\tilde{\varphi}(Z_0)\;=\;\frac{1}{4}\int_\C F\left(\frac{1}{2}(Z+Z_0)\right)\,e^{\frac{i}{2}\mu\,\mbox{\rm Im}(\overline{Z}Z_0)}\,
 \tilde{\varphi}(Z)\,{\rm d}Z\;.
$$
Applied to an element  in $\CE_\theta(\bH_1)$ of the form $F_X(Z)\;=\;e^{i\alpha \, \mbox{\rm Im}(\overline{Z}X)}$ with $X\in\C$ and $\alpha\in\R$, the above 
formula yields:
$\mbox{\rm BF}_\mu(F_X)
\tilde{\varphi}(Z_0)\;=\;e^{\frac{i\alpha}{2}\mbox{\rm Im}(\overline{Z_0}X)} \CF_\C
(\tilde{\varphi})(\mu Z_0+\alpha X)$
where $\CF_\C(\tilde{\varphi})(Z_0)\;:=\;C\int_\C e^{\frac{i}{2}\,\mbox{\rm Im}(\overline{Z}Z_0)} \tilde{\varphi}(Z)\,{\rm d}Z$.
The limit $X_0\to0$ yields $\tilde{\varphi}=\CF_\C(\tilde{\varphi})$, hence:
$$
\mbox{\rm BF}_\mu(F_X)
\tilde{\varphi}(Z_0)\;=\;e^{\frac{i\alpha}{2}\mbox{\rm Im}(\overline{Z_0}X)}\tilde{\varphi}(\mu Z_0+\alpha X)\;.
$$
\EPf
\section{Conclusions}
We now summarize what has been done in the present work. First, we establish a way to systematically produce explicit formulae for 
intertwiners of group unitary represenations. Second, applying the above intertwiner in the case of the Bargmann-Fock and Kirillov-Weyl realizations of the unitary dual of the Heisenberg groups, we realized the non-commutative torus as the range of the star-exponential for the Heisenberg group. And third, we then deduced from this a realization of the non-commutative torus as acting 
on a Fock space of entire functions.

\end{document}